\documentclass[11pt,leqno]{article}
\usepackage{amsmath, amsthm, amsfonts, amssymb}
\usepackage{mathrsfs}
\newfam\msbfam
\font\tenmsb=msbm10    \textfont\msbfam=\tenmsb \font\sevenmsb=msbm7
\scriptfont\msbfam=\sevenmsb \font\fivemsb=msbm5
\scriptscriptfont\msbfam=\fivemsb

\newfam\bigfam
\font\tenbig=msbm10 scaled \magstep2   \textfont\bigfam=\tenbig
\font\sevenbig=msbm7 scaled \magstep2 \scriptfont\bigfam=\sevenbig
\font\fivebig=msbm5 scaled \magstep2
\scriptscriptfont\bigfam=\fivebig

\def\dint{\displaystyle\int}

\def\dfrac{\displaystyle\frac}

\newtheorem{thm}{Theorem}[section]
\newtheorem{lem}{Lemma}[section]

\newtheorem{defn}{Definition}[section]
\newtheorem{pf}{Proof}

\begin{document}

{\begin{center}\bf\LARGE HERMITE-HADAMARD TYPE INEQUALITIES FOR GENERALIZED $s$-CONVEX FUNCTIONS ON REAL LINEAR FRACTAL SET  $\mathbb{R}^{\alpha}(0<\alpha<1)$

\end{center}}

\vspace{0.2cm}
 \begin{center}{ Huixia Mo\footnote[1]{
Correspondence: huixiamo@bupt.edu.cn\\
The research is supported by the NNSF of China (11471050, 11161042)} and Xin Sui

\vspace{0.2cm} School of Science, Beijing University of Posts and
Telecommunications, Beijing, 100876, China \vspace{0.2cm}\\}\end{center}

\vspace{0.3cm}

\begin{center}\begin{minipage}{5in}

\small {In the paper, we establish the Hermite-Hadamard type inequalities for the generalized s-convex functions in the second sense on  real linear fractal set $\mathbb{R}^{\alpha}(0<\alpha<1).$}\\

\vspace{0.3cm}
\end{minipage}\end{center}

\begin{center}\section{Introduction}\end{center}

The convexity of functions is an important concept in the class mathematical analysis course.  In \cite{HL}, Hudzik and Maligranda considered two kinds of convex functions which are $s$-convex. And many important inequalities are established for the $s$-convex functions. For example, the Hermite-Hadamard¡¯s inequality is one of the best known results in the literature, see \cite{KBO,AK,SYQ} and so on.

In recent years, the fractal theory has received significantly remarkable attention\cite{KST}.
The calculus on fractal set can lead to better comprehension for the various real world models from science and engineering \cite{GZS}.

On the fractal set, Mo etal.\cite{MSY,M} introduced the definition of the generalized convex function and established Hermite-Hadamard type  inequality.
In \cite{MS}, the authors introduced two kinds of generalized $s$-convex functions on fractal sets $\mathbb{R}^{\alpha}(0<\alpha<1).$

 The definitions of the generalized $s$-convex functions are as follows:

\begin{defn}\cite{MS}\;Let $\mathbb{R_+}=[0,\infty].$ A function $f:\mathbb{R_+}\rightarrow\mathbb{R^{\alpha}}$
is said to be generalized $s$-convex $(0<s<1)$ in the first sense, if
$$f(\lambda_1u+\lambda_2v)\leq\lambda_1^{\alpha s}f(u)+\lambda_2^{\alpha s}f(v), \eqno (1.1)$$
for all $u,v\in\mathbb{R_+}$ and all $\lambda_1,\lambda_2\geq0$ with $\lambda_1^ s+\lambda_2^s=1.$ We denote this by $f\in GK_s^1.$
\end{defn}

\begin{defn}\cite{MS}\;A function $f:\mathbb{R_+}\rightarrow\mathbb{R^{\alpha}}$
is said to be generalized $s$-convex $(0<s<1)$ in the second sense, if
$$f(\lambda_1u+\lambda_2v)\leq\lambda_1^{\alpha s}f(u)+\lambda_2^{\alpha s}f(v), \eqno (1.2)$$
for all $u,v\in\mathbb{R_+}$ and all $\lambda_1,\lambda_2\geq0$ with $\lambda_1 +\lambda_2=1.$ We denote this by $f\in GK_s^2.$
\end{defn}

Note that, when $s=1$, the generalized $s$-convex function in both sense is generalized  convex function\cite{MS}.

Inspired by \cite{KBO,AK,M}, in the paper we will establish the Hermite-Hadamard type inequalities for generalized $s$-convex functions.

\section{Preliminaries}

Let us review the operations with real line number on fractal space  and
use the Gao-Yang- Kang's idea to describe the definitions of the
 local fractional derivative and local fractional integral \cite{Y,WSZW,YBM,GYK,YG1}.

If $a^\alpha,b^\alpha$ and $c^\alpha$ belong to the set $\mathbb{R}^\alpha$ of real line numbers , then\\
(1)\;\;$a^\alpha+b^\alpha$ and $a^\alpha b^\alpha$ belong to the set $\mathbb{R}^\alpha$;\\
(2)\;\;$a^\alpha+b^\alpha=b^\alpha+a^\alpha=(a+b)^\alpha=(b+a)^\alpha$;\\
(3)\;\;$a^\alpha+(b^\alpha+c^\alpha)=(a^\alpha+b^\alpha)+c^\alpha$;\\
(4)\;\;$a^\alpha b^\alpha=b^\alpha a^\alpha=(ab)^\alpha=(ba)^\alpha$;\\
(5)\;\;$a^\alpha(b^\alpha c^\alpha)=(a^\alpha b^\alpha)c^\alpha$;\\
(6)\;\;$a^\alpha(b^\alpha+c^\alpha)=a^\alpha b^\alpha+a^\alpha c^\alpha$;\\
(7)\;\;$a^\alpha+0^\alpha=0^\alpha+a^\alpha=a^\alpha,$ and $a^\alpha\cdot1^\alpha=1^\alpha\cdot a^\alpha=a^\alpha$.

\begin{defn}\cite{Y}\;\;Let $f:[a,b]\rightarrow \mathbb{R^\alpha}$ be a function such that
$$ |f(x)-f(y)|<c|x-y|^\alpha \;\;\;(x,y\in [a,b]),$$
for positive constants $c$ and $\alpha(0<\alpha\leq 1).$ This function is called H\"{o}lder continuous function. In this case, we think that $f(x)$ is  in the sapce $C_\alpha[a,b].$
\end{defn}

\begin{defn}\cite{Y} The local fractional derivative of $f(x)$ of order $\alpha$ at $x=x_0$ is defined by
$$f^{(\alpha)}(x_0)=\frac{d^\alpha f(x)}{dx^\alpha}|_{x=x_0}=\lim\limits_{x\rightarrow x_0}\frac{\triangle^\alpha(f(x)-f(x_0))}{(x-x_0)^\alpha},$$
where $\triangle^\alpha(f(x)-f(x_0))\cong\Gamma(1+a)(f(x)-f(x_0)).$
If there exists $f^{((k+1)\alpha)}(x)=\overbrace{D_x^{\alpha} \ldots D_x^\alpha}\limits^{k+1\; \mbox {times}}f(x)$ for any $x\in I\subseteq\mathbb{R},$  then we denoted
$f\in D_{(k+1)\alpha}(I)$, where $k=0,1,2\dots.$
\end{defn}

\begin{defn}\cite{Y}\;\;Let $f\in C_{\alpha}[a,b]$. Then the local fractional integral of the function $f(x)$ is defined by,
$$\begin{array}{cl}&_aI_b^{(\alpha)}f(x)\\
&
=\dfrac{1}{\Gamma(1+a)}\int_a^bf(t)(dt)^\alpha\\&
=\dfrac{1}{\Gamma(1+a)}\lim\limits_{\triangle t\rightarrow 0}\sum\limits_{j=0}^{N}f(t_j)(\bigtriangleup t_j)^\alpha,\end{array}$$
where $\triangle t_j=t_{j+1}-t_j,\;\triangle t=\max\{\triangle t_1,\triangle t_2,\triangle t_j,\ldots\}$ and $[t_j,t_j+1],$ $j=0,\ldots,N-1,$ $t_0=a,$ $t_N=b,$ is a partition of the interval $[a,b]$.
\end{defn}

\begin{lem}\cite{Y}\;\;Suppose that $g\in C_1[a,b]$ and  $f\in C_\alpha[g(a), g(b)],$ then
$$_{g(a)}I_{g(b)}^{(\alpha)}f(x)=_{a}I_{b}^{(\alpha)}f(g)(s)[g'(s)]^{\alpha}.$$
 \end{lem}

\begin{lem}\cite{Y}\\
(1)\;\; Suppose that $f(x)=g^{(\alpha)}(x)\in C_{\alpha}[a,b]$, then we have
   $$_{a}I_{b}^{(\alpha)}f(x)=g(b)-g(a).$$
(2)\;\;
 Suppose that $f(x), g(x)\in D_{\alpha}[a,b]$ and $f^{(\alpha)}(x), g^{(\alpha)}(x)\in C_{\alpha}[a,b]$, then we have
   $$_{a}I_{b}^{\alpha}f(x)g^{(\alpha)}(x)=f(x)g(x)\bigg|_{a}^{b}-_aI_{b}^{(\alpha)}f^{(\alpha)}(x)g(x).$$
 \end{lem}

\begin{lem}\cite{Y}\;\;\\
$$\dfrac{d^\alpha x^{k\alpha}}{dx^\alpha}=\dfrac{\Gamma(1+k\alpha)}{\Gamma(1+(k-1)\alpha)}x^{(k-1)\alpha}.$$

From the above formula and Lemma 2.2, we have
$$\dfrac{1}{\Gamma(1+\alpha)}\int_a^bx^{k\alpha}(dx)^{\alpha}=\dfrac{\Gamma(1+k\alpha)}{\Gamma(1+(k+1)\alpha)}(b^{(k+1)\alpha}-a^{(k+1)\alpha}),k\in R.$$
\end{lem}

\begin{lem}\cite{Y}\;\; \textbf{(Generalized H\"{o}lder¡¯s inequality)}
 Let $f,g\in C_{\alpha}[a,b]$, $p,q>1$ with $1/p+1/q=1$, then
 $$\dfrac{1}{\Gamma(1+\alpha)}\int_a^b|f(x)g(x)|(dx)^{\alpha}\leq\bigg(\dfrac{1}{\Gamma(1+\alpha)}\int_a^b|f(x)|^p(dx)^{\alpha}\bigg)^{1/p}\bigg(\dfrac{1}{\Gamma(1+\alpha)}\int_a^b|g(x)|^q(dx)^{\alpha}\bigg)^{1/q}.$$
\end{lem}

\begin{center}\section{Main Results}\end{center}

\begin{thm} Suppose that $f:\mathbb{R_+}\rightarrow\mathbb{R^{\alpha}}$ is a generalized $s$-convex function in the second sense, where $s\in(0,1).$ Let $a,b\in[0,\infty),$ $a<b.$ If $f\in C_{\alpha}[a,b]$, then the following inequalities hold:
$$\begin{array}{cl}\dfrac{2^{(s-1)\alpha}}{\Gamma(1+\alpha)}f\biggl(\dfrac{a+b}{2}\biggr)\leq
\dfrac{_{a}I_{b}^{\alpha}f(x)}{(b-a)^\alpha}\leq\dfrac{\Gamma(1+s\alpha)}{\Gamma(1+(s+1)\alpha)}\bigl(f(a)+f(b)\bigr).\end{array}\eqno{(3.1)}$$
\end{thm}

\begin{pf}
Let $x=a+b-t,$ then by Lemma 2.1, we have $_{\frac{a+b}{2}}I_b^{(\alpha)}f(x)={_a}I_{\frac{a+b}{2}}^{(\alpha)}f(a+b-t).$

 Moreover, $f$ is $s$ convex in the second sense, then
$$\begin{array}{cl}_aI_b^{(\alpha)}f(x)
&
={_a}I_{\frac{a+b}{2}}^{(\alpha)}\bigl[f(x)+f(a+b-x)\bigr]\\&
\geq 2^{\alpha s} {_a}I_{\frac{a+b}{2}}^{(\alpha)}f\biggl(\dfrac{a+b}{2}\biggr)\\&
=\dfrac{2^{(s-1)\alpha}}{\Gamma(1+\alpha)}(b-a)^\alpha f\biggl(\dfrac{a+b}{2}\biggr).\end{array}$$

In the other hand, let $x=b-(b-a)t$, $0\leq t\leq 1,$ then we get
$$\begin{array}{cl}_aI_b^{(\alpha)}f(x)&=(b-a)^\alpha {_0}I_1^{(\alpha)}f\bigl[ta+(1-t)b\bigr]\\&
\leq(b-a)^\alpha {_0}I_1^{(\alpha)}\bigl[t^{\alpha s}f(a)+(1-t)^{\alpha s}f(b)\bigr]\\&
=(b-a)\bigl[f(a) {_0}I_1^{(\alpha)}t^{\alpha s}+f(b) {_0}I_1^{(\alpha)}(1-t)^{\alpha s}\bigr].\end{array}$$

From Lemma 2.3, it is easy to see that
$${_0}I_1^{(\alpha)}t^{\alpha s}=\dfrac{\Gamma(1+s\alpha)}{\Gamma(1+(s+1)\alpha)},$$
 and
 $${_0}I_1^{(\alpha)}(1-t)^{\alpha s}=\dfrac{\Gamma(1+s\alpha)}{\Gamma(1+(s+1)\alpha)}.$$

 So,
$$_aI_b^{(\alpha)}f(x)\leq(b-a)^\alpha\dfrac{\Gamma(1+s\alpha)}{\Gamma(1+(s+1)\alpha)}\bigl(f(a)+f(b)\bigr).$$

Therefore, we obtain
$$\dfrac{2^{(s-1)\alpha}}{\Gamma(1+\alpha)}f\biggl(\dfrac{a+b}{2}\biggr)\leq
\dfrac{_{a}I_{b}^{\alpha}f(x)}{(b-a)^\alpha}\leq\dfrac{\Gamma(1+s\alpha)}{\Gamma(1+(s+1)\alpha)}\bigl(f(a)+f(b)\bigr).$$
\end{pf}

\begin{thm}
Let $f: I\rightarrow \mathbb{R^{\alpha}}$, be a differentiable function on $I^0$($I^0$ is the interior of $I$) such that $f^{(\alpha)}\in C_{\alpha}[a,b]$, where $a,b\in I^0$, $a<b$. If $|f^{(\alpha)}|^q$ is generalized $s$-convex in the second sense on $[a,b]$ for some fixed $s\in(0,1)$ and $q\geq 1$, then
$$\begin{array}{cl}&\bigg|\dfrac{f(a)+f(b)}{2^\alpha}-\dfrac{\Gamma(1+\alpha)}{(b-a)^\alpha}{_a}I_b^{(\alpha)}f(x)\bigg|\\
&\leq\dfrac{(b-a)^\alpha}{2^\alpha}\bigg(\dfrac{\Gamma(1+\alpha)}{\Gamma(1+2\alpha)}\bigg)^{\frac{q-1}{q}}\bigg[\dfrac{\Gamma(1+s\alpha)}{\Gamma(1+(s+1)\alpha)}+\dfrac{\Gamma(1+\alpha)\Gamma(1+s\alpha)}{\Gamma(1+(s+2)\alpha)}
\biggl(\biggl(\dfrac{1}{2}\biggr)^{\alpha s}-2^\alpha\biggr)\bigg]^{\frac{1}{q}}\\
&\quad\times\bigg[\big|f^{(\alpha)}(a)\big|^q+\big|f^{(\alpha)}(b)\big|^q\bigg]^{\frac{1}{q}}\\
\end{array}$$

\end{thm}

To show Theorem 3.2  is right, we need the following Lemma.

\begin{lem}(\cite{M})\;\; Let $I\subseteq R$ be an interval, $f:I^{0}\subseteq R\rightarrow R^{\alpha}$  such that $f\in D_\alpha(I^0)$ and $f^{(\alpha)}\in C_\alpha[a,b]$ for $a,b\in I^0$ with $a<b$. Then the following equality holds:
$$\dfrac{f(a)+f(b)}{2^\alpha}-\dfrac{\Gamma(1+\alpha)}{(b-a)^\alpha}{_a}I_b^{(\alpha)}f(x)=\dfrac{(b-a)^\alpha}{2^\alpha}\dfrac{1}{\Gamma(1+\alpha)}\int_0^1(1-2t)^\alpha f^{(\alpha)}\big(ta+(1-t)b\big)(dt)^\alpha.$$
\end{lem}

Now, let us give the proof of Theorem 3.2.
\begin{pf}
 From Lemma 3.1, it is obvious that
$$\begin{array}{cl}
&\bigg|\dfrac{f(a)+f(b)}{2^\alpha}-\dfrac{\Gamma(1+\alpha)}{(b-a)^\alpha}{_a}I_b^{(\alpha)}f(x)\bigg|\\
\leq&\dfrac{(b-a)^\alpha}{2^\alpha}\dfrac{1}{\Gamma(1+\alpha)}\int_0^1|1-2t|^\alpha \big|f^{(\alpha)}(ta+(1-t)b)\big|(dt)^\alpha.\end{array}\eqno{(3.2)}$$

Let us estimate $$\dfrac{1}{\Gamma(1+\alpha)}\int_0^1|1-2t|^\alpha \big|f^{(\alpha)}(ta+(1-t)b)\big|(dt)^\alpha,$$ for $q=1$ and $q>1.$

{\bf Case I.}\;\; $q=1.$

 Since $|f^{(\alpha)}|$ is generalized $s$-convex on $[a,b]$ in the second sense, we can know that for any $t\in [0,1]$
$$\big|f^{(\alpha)}(ta+(1-t)b)\big|\leq t^{\alpha s}|f^{(\alpha)}(a)|+(1-t)^{\alpha s}|f^{(\alpha)}(b)|.$$

Then,  we have,
$$\begin{array}{cl}&\dfrac{1}{\Gamma(1+\alpha)}\int_0^1|1-2t|^\alpha \big|f^{(\alpha)}(ta+(1-t)b)\big|(dt)^\alpha\\
&\leq\dfrac{1}{\Gamma(1+\alpha)}\int_0^1|1-2t|^\alpha \big[t^{\alpha s}|f^{(\alpha)}(a)|+(1-t)^{\alpha s}|f^{(\alpha)}(b)|\big](dt)^\alpha\\
&=\biggl[|f^{(\alpha)}(a)|\dfrac{1}{\Gamma(1+\alpha)}\int_0^1t^{\alpha s}|1-2t|^\alpha(dt)^\alpha+|f^{(\alpha)}(b)|\dfrac{1}{\Gamma(1+\alpha)}\int_0^1(1-t)^{\alpha s}|1-2t|^{\alpha}(dt)^\alpha\biggr].\\
\end{array}\eqno{(3.3)}$$

From Lemma 2.2 and lemma 2.3, it is easy to see that
$$\begin{array}{cl}
&\dfrac{1}{\Gamma(1+\alpha)}\int_0^1t^{\alpha s}|1-2t|^\alpha(dt)^\alpha\\
&=\dfrac{1}{\Gamma(1+\alpha)}\bigg[\int_0^{\frac{1}{2}}t^{\alpha s}(1-2t)^\alpha(dt)^\alpha
+\dint_{\frac{1}{2}}^1t^{\alpha s}(2t-1)^\alpha(dt)^\alpha\bigg]\\
&=\bigg[\dfrac{\Gamma(1+s\alpha)}{\Gamma(1+(s+1)\alpha)}+\dfrac{\Gamma(1+\alpha)\Gamma(1+s\alpha)}{\Gamma(1+(s+2)\alpha)}\biggl(\bigg(\dfrac{1}{2}\biggr)^{\alpha s}-2^\alpha\biggr)\bigg].
\end{array}\eqno{(3.4)}$$
And, let $1-t=x,$ then by Lemma 2.1 and (3.4), we have
$$\begin{array}{cl}
&\dfrac{1}{\Gamma(1+\alpha)}\int_0^1(1-t)^{\alpha s}|1-2t|^\alpha(dt)^\alpha\\
&=\dfrac{1}{\Gamma(1+\alpha)}\int_0^1 x^{\alpha s}|1-2x|^\alpha(dx)^\alpha\\
&=\bigg[\dfrac{\Gamma(1+s\alpha)}{\Gamma(1+(s+1)\alpha)}+\dfrac{\Gamma(1+\alpha)\Gamma(1+s\alpha)}{\Gamma(1+(s+2)\alpha)}\biggl(\bigg(\dfrac{1}{2}\biggr)^{\alpha s}-2^\alpha\biggr)\bigg].
\end{array}\eqno{(3.5)}$$

Thus, substituting (3.4) and (3.5) into (3.3), we have

$$\begin{array}{cl}&\dfrac{1}{\Gamma(1+\alpha)}\int_0^1|1-2t|^\alpha \big|f^{(\alpha)}(ta+(1-t)b)\big|(dt)^\alpha\\
&\leq\dfrac{\Gamma(1+\alpha)}{\Gamma(1+2\alpha)}
\bigg[\dfrac{\Gamma(1+s\alpha)}{\Gamma(1+(s+1)\alpha)}+\dfrac{\Gamma(1+\alpha)\Gamma(1+s\alpha)}{\Gamma(1+(s+2)\alpha)}
\biggl(\biggl(\dfrac{1}{2}\biggr)^{\alpha s}-2^\alpha\biggr)\bigg]\\
&\quad\times\bigg[\big|f^{(\alpha)}(a)\big|^q+\big|f^{(\alpha)}(b)\big|\bigg].\\
\end{array}\eqno{(3.6)}$$

Thus, from (3.2), we obtain

$$\begin{array}{cl}&\bigg|\dfrac{f(a)+f(b)}{2^\alpha}-\dfrac{\Gamma(1+\alpha)}{(b-a)^\alpha}{_a}I_b^{(\alpha)}f(x)\bigg|\\
&\leq\dfrac{(b-a)^\alpha}{2^\alpha}\dfrac{\Gamma(1+\alpha)}{\Gamma(1+2\alpha)}
\bigg[\dfrac{\Gamma(1+s\alpha)}{\Gamma(1+(s+1)\alpha)}+\dfrac{\Gamma(1+\alpha)\Gamma(1+s\alpha)}{\Gamma(1+(s+2)\alpha)}
\biggl(\biggl(\dfrac{1}{2}\biggr)^{\alpha s}-2^\alpha\biggr)\bigg]\\
&\quad\times\bigg[\big|f^{(\alpha)}(a)\big|^q+\big|f^{(\alpha)}(b)\big|\bigg].\\
\end{array}$$

{\bf Case II.}\;\; $q>1.$

Using the generalized H\"{o}lder's inequality(Lemma 2.4), we obtain
$$\begin{array}{cl}&\dfrac{1}{\Gamma(1+\alpha)}\int_0^1|1-2t|^\alpha\big|f^{(\alpha)}(ta+(1-t)b)\big|(dt)^\alpha\\&
=\dfrac{1}{\Gamma(1+\alpha)}\int_0^1|1-2t|^{\alpha\frac{q-1}{q}}|1-2t|^{\alpha\frac{1}{q}}\big|f^{(\alpha)}(ta+(1-t)b)\big|(dt)^\alpha\\
&\leq\bigg(\dfrac{1}{\Gamma(1+\alpha)}\int_0^1|1-2t|^\alpha(dt)^\alpha\bigg)^{\frac{q-1}{q}}\bigg(\dfrac{1}
{\Gamma(1+\alpha)}\int_0^1|1-2t|^\alpha|f^{(\alpha)}(ta+(1-t)b)|^q(dt)^\alpha\bigg)^{\frac{1}{q}}.\end{array}\eqno{(3.7)}$$

It is obvious that
$$\begin{array}{cl}&\dfrac{1}{\Gamma(1+\alpha)}\int_0^1|1-2t|^\alpha(dt)^\alpha\\
&=\dfrac{1}{\Gamma(1+\alpha)}\int_0^{\frac{1}{2}}(1-2t)^\alpha(dt)^\alpha+
\dfrac{1}{\Gamma(1+\alpha)}\int_{\frac{1}{2}}^1(2t-1)^\alpha(dt)^\alpha=\dfrac{\Gamma(1+\alpha)}{\Gamma(1+2\alpha)}.\end{array}\eqno{(3.8)}$$

 Moreover, since $|f^{(\alpha)}|^q$ is generalized $s$-convex in the second sense on $[a,b],$ then
 $$\begin{array}{cl}&\dfrac{1}{\Gamma(1+\alpha)}\int_0^1|1-2t|^\alpha\big|f^{(\alpha)}(ta+(1-t)b)\big|^q(dt)^\alpha\\
 &\leq\dfrac{1}{\Gamma(1+\alpha)}\int_0^1|1-2t|^\alpha\bigg(t^{\alpha s}|f^{(\alpha)}(a)|^q+(1-t)^{\alpha s}|f^{(\alpha)}(b)|^q\bigg)(dt)^\alpha\\
&=|f^{(\alpha)}(a)|^q\dfrac{1}{\Gamma(1+\alpha)}\int_0^1|1-2t|^\alpha t^{\alpha s}(dt)^\alpha+|f^{(\alpha)}(b)|^q\dfrac{1}{\Gamma(1+\alpha)}\int_0^1|1-2t|^\alpha(1-t)^{\alpha s}(dt)^\alpha.\\
\end{array}$$

  From (3.3) and (3.4), it is easy to see that

$$\begin{array}{cl}
&\dfrac{1}{\Gamma(1+\alpha)}\int_0^1|1-2t|^\alpha t^{\alpha s}(dt)^\alpha\\
=&\dfrac{1}{\Gamma(1+\alpha)}\int_0^1|1-2t|^\alpha(1-t)^{\alpha s}(dt)^\alpha\\
=&\bigg[\dfrac{\Gamma(1+s\alpha)}{\Gamma(1+(s+1)\alpha)}+\dfrac{\Gamma(1+\alpha)\Gamma(1+s\alpha)}{\Gamma(1+(s+2)\alpha)}\biggl(\biggl(\dfrac{1}{2}\biggr)^{\alpha s}-2^\alpha\biggr)\bigg].
\end{array}$$

So,  $$\begin{array}{cl}&\dfrac{1}{\Gamma(1+\alpha)}\int_0^1|1-2t|^\alpha\big|f^{(\alpha)}(ta+(1-t)b)\big|^q(dt)^\alpha\\
 =&\bigg[\dfrac{\Gamma(1+s\alpha)}{\Gamma(1+(s+1)\alpha)}+\dfrac{\Gamma(1+\alpha)\Gamma(1+s\alpha)}{\Gamma(1+(s+2)\alpha)}
\biggl(\biggl(\dfrac{1}{2}\biggr)^{\alpha s}-2^\alpha\biggr)\bigg]^{\frac{1}{q}}\bigg[\big|f^{(\alpha)}(a)\big|^q+\big|f^{(\alpha)}(b)\big|^q\bigg]^{\frac{1}{q}}.
\end{array}\eqno{(3.9)}$$

Thus, substituting (3.8) and (3.9) into (3.7), we have

$$\begin{array}{cl}&\dfrac{1}{\Gamma(1+\alpha)}\int_0^1|1-2t|^\alpha\big|f^{(\alpha)}(ta+(1-t)b)\big|(dt)^\alpha\\
\leq&\bigg(\dfrac{\Gamma(1+\alpha)}{\Gamma(1+2\alpha)}\bigg)^{\frac{q-1}{q}}\bigg[\dfrac{\Gamma(1+s\alpha)}{\Gamma(1+(s+1)\alpha)}+\dfrac{\Gamma(1+\alpha)\Gamma(1+s\alpha)}{\Gamma(1+(s+2)\alpha)}
\biggl(\biggl(\dfrac{1}{2}\biggr)^{\alpha s}-2^\alpha\biggr)\bigg]^{\frac{1}{q}}
\bigg[\big|f^{(\alpha)}(a)\big|^q+\big|f^{(\alpha)}(b)\big|^q\bigg]^{\frac{1}{q}}.\end{array}$$

So, from (3.2) it follows that
$$\begin{array}{cl}&\bigg|\dfrac{f(a)+f(b)}{2^\alpha}-\dfrac{\Gamma(1+\alpha)}{(b-a)^\alpha}{_a}I_b^{(\alpha)}f(x)\bigg|\\
&\leq\dfrac{(b-a)^\alpha}{2^\alpha}\bigg(\dfrac{\Gamma(1+\alpha)}{\Gamma(1+2\alpha)}\bigg)^{\frac{q-1}{q}}\bigg[\dfrac{\Gamma(1+s\alpha)}{\Gamma(1+(s+1)\alpha)}+\dfrac{\Gamma(1+\alpha)\Gamma(1+s\alpha)}{\Gamma(1+(s+2)\alpha)}
\biggl(\biggl(\dfrac{1}{2}\biggr)^{\alpha s}-2^\alpha\biggr)\bigg]^{\frac{1}{q}}\\
&\quad\times\bigg[\big|f^{(\alpha)}(a)\big|^q+\big|f^{(\alpha)}(b)\big|^q\bigg]^{\frac{1}{q}}.\\
\end{array}$$

Thus, we complete the proof of Theorem 3.2.
\end{pf}

\begin{thm}Let $f:$ $I\rightarrow \mathbb{R^{\alpha}},$ $I\subset[0,\infty)$, be a differentiable function on $I^0$ such that $f^{(\alpha)}\in C_{\alpha}[a,b]$, where $a,b\in I$ and $a<b$. If $|f^{(\alpha)}|^q$ is generalized $s$-convex in the second sense on $[a,b]$ for some fixed $s\in(0,1)$ and $q>1$, then
$$\begin{array}{cl}&\bigg|\dfrac{f(a)+f(b)}{2^\alpha}-\dfrac{\Gamma(1+\alpha)}{(b-a)^\alpha}{_a}I_b^{(\alpha)}f(x)\bigg|\\
&\leq\dfrac{(b-a)^\alpha}{2^\alpha}\bigg[\dfrac{\Gamma(1+\frac{q}{q-1}\alpha)}{2^\alpha\Gamma(1+\frac{2q-1}{q-1}\alpha)}\bigg]^{\frac{q-1}{q}}
\bigg(\dfrac{\Gamma(1+s\alpha)}{2^\alpha\Gamma(1+(s+1)\alpha)}\bigg)^{\frac{1}{q}}\\&
\quad\times\biggl[\biggl(|f^{(\alpha)}(a)|^q+\biggl|f^{(\alpha)}\biggl(\dfrac{a+b}{2}\biggr)\biggr|^q\biggr)^{\frac{1}{q}}+\biggl(\biggl|f^{(\alpha)}\biggl(\dfrac{a+b}{2}\biggr)\biggr|^q+|f^{(\alpha)}(b)|^q\biggr)^{\frac{1}{q}}\biggr].
\end{array}$$
\end{thm}

\begin{pf}
From Lemma 3.1, we have
$$\begin{array}{cl}&\bigg|\dfrac{f(a)+f(b)}{2^\alpha}-\dfrac{\Gamma(1+\alpha)}{(b-a)^\alpha}{_a}I_b^{(\alpha)}f(x)\bigg|\\
&\leq\dfrac{(b-a)^\alpha}{2^\alpha}\dfrac{1}{\Gamma(1+\alpha)}\int_0^1|1-2t|^\alpha |f^{(\alpha)}(ta+(1-t)b)|(dt)^\alpha\\
&\leq\dfrac{(b-a)^\alpha}{2^\alpha}\bigg[\dfrac{1}{\Gamma(1+\alpha)}\int_0^{\frac{1}{2}}(1-2t)|f^{(\alpha)}(ta+(1-t)b)|(dt)^\alpha\\
&\quad+\dfrac{1}{\Gamma(1+\alpha)}\int_{\frac{1}{2}}^1(2t-1)|f^{(\alpha)}(ta+(1-t)b)|(dt)^\alpha\bigg]
.\end{array}\eqno{(3.10)}$$

Let us estimate
$$\dfrac{1}{\Gamma(1+\alpha)}\int_0^{\frac{1}{2}}(1-2t)|f^{(\alpha)}(ta+(1-t)b)|(dt)^\alpha$$
and
$$\dfrac{1}{\Gamma(1+\alpha)}\int_{\frac{1}{2}}^1(2t-1)|f^{(\alpha)}(ta+(1-t)b)|(dt)^\alpha,$$
respectively.

 Using the generalized H\"{o}lder's inequality(Lemma 2.4), we obtain
$$\begin{array}{cl}&\dfrac{1}{\Gamma(1+\alpha)}\int_0^{\frac{1}{2}}(1-2t)|f^{(\alpha)}(ta+(1-t)b)|(dt)^\alpha\\
&\leq\bigg(\dfrac{1}{\Gamma(1+\alpha)}\int_0^{\frac{1}{2}}(1-2t)^{\alpha\frac{q}{q-1}}(dt)^\alpha\bigg)^{\frac{q-1}{q}}
\bigg(\dfrac{1}{\Gamma(1+\alpha)}\int_0^{\frac{1}{2}}|f^{(\alpha)}(ta+(1-t)b)|^{\alpha q}(dt)^\alpha\bigg)^{\frac{1}{q}}.\end{array}\eqno{(3.11)}$$

It is easy to see that
$$\begin{array}{cl}\dfrac{1}{\Gamma(1+\alpha)}\int_0^{\frac{1}{2}}(1-2t)^{\alpha\frac{q}{q-1}}(dt)^\alpha=\dfrac{1}{\Gamma(1+\alpha)}\int_{\frac{1}{2}}^1(1-2t)^{\alpha\frac{q}{q-1}}(dt)^\alpha
=\frac{\Gamma(1+\frac{q}{q-1}\alpha)}{2^\alpha\Gamma(1+\frac{2q-1}{q-1}\alpha)}.\end{array}\eqno{(3.12)}$$

Let $|f^{(\alpha)}(ta+(1-t)b)|^q=U(t).$ Since $|f^{(\alpha)}|^q$ is a generalized $s$ convex function in the second sense , then $U(t)$ is  a generalized $s$ convex function in the second sense , too. Thus, from the right hand side of (3.1),
 it follows that

$$\begin{array}{cl}&\dfrac{1}{\Gamma(1+\alpha)}\int_0^{\frac{1}{2}}|f^{(\alpha)}(ta+(1-t)b)|^q(dt)^\alpha \\ =&\dfrac{1}{\Gamma(1+\alpha)}\int_0^{\frac{1}{2}}U(t)(dt)^\alpha\\
\leq&\biggl(\dfrac{1}{2}-0\biggr)^\alpha\dfrac{\Gamma(1+s\alpha)}{\Gamma(1+(s+1)\alpha)}(U(0)+U(\frac{1}{2}))\\
=&\dfrac{\Gamma(1+s\alpha)}
  {2^\alpha\Gamma(1+(s+1)\alpha)}\biggl(\biggl|f^{(\alpha)}\biggl(\dfrac{a+b}{2}\biggr)\biggr|^q+|f^{(\alpha)}(b)|^q\biggr).\end{array}\eqno{(3.13)}$$

Thus, substituting (3.12) and (3.13) into (3.11), we get

$$\begin{array}{cl}&\dfrac{1}{\Gamma(1+\alpha)}\int_0^{\frac{1}{2}}(1-2t)|f^{(\alpha)}(ta+(1-t)b)|(dt)^\alpha\\
&\leq\bigg(\dfrac{\Gamma(1+\frac{q}{q-1}\alpha)}{2^\alpha\Gamma(1+\frac{2q-1}{q-1}\alpha)}\bigg)^{\frac{q-1}{q}}
\bigg(\dfrac{\Gamma(1+s\alpha)}{2^\alpha\Gamma(1+(s+1)\alpha)}\bigg)^{\frac{1}{q}}
\biggl(|f^{(\alpha)}(a)|^q+\biggl|f^{(\alpha)}\biggl(\frac{a+b}{2}\biggr)\biggr|^q\biggr)^{\frac{1}{q}}.\end{array}\eqno{(3.14)}$$

Moreover,
$$\begin{array}{cl}\dfrac{1}{\Gamma(1+\alpha)}\int_{\frac{1}{2}}^1(2t-1)^{\alpha\frac{q}{q-1}}(dt)^\alpha=\dfrac{1}{\Gamma(1+\alpha)}\int_{\frac{1}{2}}^1(1-2t)^{\alpha\frac{q}{q-1}}(dt)^\alpha
=\frac{\Gamma(1+\frac{q}{q-1}\alpha)}{2^\alpha\Gamma(1+\frac{2q-1}{q-1}\alpha)}.\end{array}$$

And, similar to the estimate of (3.13), we have
   $$\begin{array}{cl}\dfrac{1}{\Gamma(1+\alpha)}\int_{\frac{1}{2}}^1|f^{(\alpha)}(ta+(1-t)b)|^q(dt)^\alpha\leq\dfrac{\Gamma(1+s\alpha)}
  {2^\alpha\Gamma(1+(s+1)\alpha)}\biggl(\biggl|f^{(\alpha)}\biggl(\dfrac{a+b}{2}\biggr)\biggr|^q+|f^{(\alpha)}(b)|^q\biggr).\end{array}$$

So, it is analogues to  the estimate of (3.11), we have

$$\begin{array}{cl}&\dfrac{1}{\Gamma(1+\alpha)}\int_0^{\frac{1}{2}}(1-2t)|f^{(\alpha)}(ta+(1-t)b)|(dt)^\alpha\\
&\leq\bigg(\dfrac{1}{\Gamma(1+\alpha)}\int_{\frac{1}{2}}^1(2t-1)^{\alpha\frac{q}{q-1}}(dt)^\alpha\bigg)^{\frac{q-1}{q}}
\bigg(\dfrac{1}{\Gamma(1+\alpha)}\int_{\frac{1}{2}}^1|f^{(\alpha)}(ta+(1-t)b)|^{q\alpha}(dt)^\alpha\bigg)^{\frac{1}{q}}\\
&\leq\bigg(\dfrac{\Gamma(1+\frac{q}{q-1}\alpha)}{2^\alpha\Gamma(1+\frac{2q-1}{q-1}\alpha)}\bigg)^{\frac{q-1}{q}}
\bigg(\dfrac{\Gamma(1+s\alpha)}{2^\alpha\Gamma(1+(s+1)\alpha)}\bigg)^{\frac{1}{q}}
\biggl(\biggl|f^{(\alpha)}\biggl(\frac{a+b}{2}\biggr)\biggr|^q+|f^{(\alpha)}(b)|^q\biggr)\biggr|^q\biggr)^{\frac{1}{q}}.\end{array}\eqno{(3.15)}$$

Thus, combining (3.10), (3.14) and (3.15), we obtain
$$\begin{array}{cl}&\bigg|\dfrac{f(a)+f(b)}{2^\alpha}-\dfrac{\Gamma(1+\alpha)}{(b-a)^\alpha}{_a}I_b^{(\alpha)}f(x)\bigg|\\
&\leq\dfrac{(b-a)^\alpha}{2^\alpha}\bigg[\dfrac{\Gamma(1+\frac{q}{q-1}\alpha)}{2^\alpha\Gamma(1+\frac{2q-1}{q-1}\alpha)}\bigg]^{\frac{q-1}{q}}
\bigg(\dfrac{\Gamma(1+s\alpha)}{2^\alpha\Gamma(1+(s+1)\alpha)}\bigg)^{\frac{1}{q}}\\&
\quad\times\biggl[\biggl(|f^{(\alpha)}(a)|^q+\biggl|f^{(\alpha)}\biggl(\dfrac{a+b}{2}\biggr)\biggr|^q\biggr)^{\frac{1}{q}}+\biggl(\biggl|f^{(\alpha)}\biggl(\dfrac{a+b}{2}\biggr)\biggr|^q+|f^{(\alpha)}(b)|^q\biggr)^{\frac{1}{q}}\biggr].
\end{array}$$
\end{pf}

Terefore, we complete the proof of Theorem 3.3.

\noindent{\bf Acknowledgments}\;\;{This work is supported by the National Natural Science Foundation of China (11161042, 11471050)}

\end{document}